\documentclass[12pt]{amsart}
\usepackage[margin=1.3in]{geometry}
\usepackage{url}
\usepackage[parfill]{parskip}
\usepackage{graphicx}
\usepackage{amsmath}
\usepackage{amssymb}
\usepackage{epstopdf}
\title[Bounds on the number of integer points in a polytope]
{Bounds on the number of integer points in a polytope via concentration estimates}
\author{Austin Shapiro}
\date{{\small October 2010}}
\address{{\small Department of Mathematics, University of Michigan, Ann Arbor, MI 48109-1043, USA}}
\thanks{{\small The author was supported in part by NSF grants DMS 0400617 and DMS 0856640.}}

\begin{document}
\maketitle
\newtheorem{Lem}{Lemma}
\newtheorem{Cor}{Corollary}
\newtheorem{Cor'}{Corollary}
\newtheorem{Thm}{Theorem}
\newtheorem{Thm'}{Theorem}
\newtheorem{Prop}{Proposition}
\newtheorem{Prop'}{Proposition}
\newtheorem{Def}{Definition}
\newtheorem{Defs}[Def]{Definitions}
\newtheorem*{Thm2a}{Theorem 2a}
\newtheorem*{Thm2b}{Theorem 2b}
\newtheorem*{Prop3a}{Proposition 3a}

\begin{abstract}
It is generally hard to count, or even estimate, how many integer points lie in a polytope $P$. Barvinok and Hartigan have approached the problem by way of information theory, showing how to efficiently compute a random vector which samples the integer points of $P$ with (computable) constant mass, but which may also land outside $P$. Thus, to count the integer points of $P$, it suffices to determine the frequency with which the random vector falls in $P$.

We prove a collection of efficiently computable upper bounds on this frequency. We also show that if $P$ is suitably presented by $n$ linear inequalities and $m$ linear equations ($m$ fixed), then under mild conditions separating the expected value of the above random vector from the origin, the frequency with which it falls in $P$ is $O(n^{-m/2})$ as $n\rightarrow\infty$. As in the classical {\sl Littlewood--Offord problem}, all results in the paper are obtained by bounding the point concentration of a sum of independent random variables; we sketch connections to previous work on the subject.
\end{abstract}

\section{Introduction}
The problem of counting integer points in polytopes has been extensively studied, and appears to be quite difficult in general. It is NP-hard to determine whether an arbitrary integral polytope with $n$ facets contains an integer point at all~\cite{GJ1979}. Given this state of affairs, attention has largely shifted to approximating or bounding the number of integer points in a polytope, and the closely related problem of sampling almost uniformly from the set of integer points in a polytope.

For certain classes of polytopes,
almost uniform sampling has been achieved by specially constructed Markov chains with good mixing properties. One notable success of this method is due to Jerrum, Sinclair, and Vigoda, who in~\cite{JSV2004} construct a fully polynomial randomized approximation scheme for the permanent of a 0-1 matrix (equal to the number of integer points in a perfect matching polytope). However, for a general polytope $P$, it is not known how to efficiently generate Markov chains which sample almost uniformly from the integer points in $P$. A survey of this and other approaches to the problem can be found in~\cite{DeLoera2005},~\cite{DeLoera2009}.

In~\cite{BH2010}, Barvinok and Hartigan proposed a new approach to the problem using the principle of maximum entropy. Given a polytope $P\subset{\Bbb R}^n$ defined by the inequalities
$$x_1\geq 0,~ x_2 \geq 0,~ \ldots,~ x_n\geq 0, \quad A{\bf x} = {\bf b},$$
where $A\in{\Bbb R}^{m\times n}$ and ${\bf b}\in{\Bbb R}^m$, they introduce a random vector $X = (X_1,X_2,\ldots,X_n)$ of maximum entropy, subject to the constraints that all coordinates are distributed on ${\bf Z}_{\geq 0}$ and that ${\bf E}[AX] = {\bf b}$ (i.e., the mean of $X$ lies in $P$). This random vector has constant mass $e^{-H(X)}$ on all points of $P\cap {\Bbb Z}^n$, where $H(X)$ is the {\sl entropy} of the random vector, defined by
$$H(X) ~:=~ -\sum_{k_1,\ldots,k_n\geq 0} {\bf Pr}[X=(k_1,\ldots,k_n)]~ \ln {\bf Pr}[X=(k_1,\ldots,k_n)].$$
Thus $X$ is, in a sense, a good approximation of the uniform distribution on $P\cap {\Bbb Z}^n$. However, not all of the mass of $X$ lies in $P$; thus we have
$$|P\cap {\Bbb Z}^n| \quad = \quad e^{H(X)} {\bf Pr}[X\in P].$$

As it turns out~\cite{BH2010}, the coordinates of $X$ are independent and {\sl geometrically distributed}, that is, there exist $q_j\in [0,1),~~1\leq j\leq n$ so that
$${\bf Pr}[X_j=k] = (1-q_j)q_j^k \quad{\rm for}\quad k\in{\Bbb Z}_{\geq 0}.$$
After a change of parameter $z_j:={\bf E}[X_j]=\frac{q_j}{1-q_j}$, the entropy $H(X)$ may be written as
\begin{equation}
H(X) = \sum_{j=1}^n (z_j+1)\ln(z_j+1)-z_j\ln z_j.
\end{equation}
This is a strictly concave function of $z_1,\ldots,z_n$, so it can be maximized efficiently by (e.g.) interior point methods (for details, see~\cite{BH2010}). Thus the parameters $q_j$, and with them the distribution and entropy of $X$, are efficiently computable. Hence, the outstanding question is how to bound the factor ${\bf Pr}[X\in P]$, particularly under weak assumptions (i.e., when a local central limit theorem is not feasible). This paper offers several upper bounds.

\section{Summary of results}
\subsection{Definitions and notation}
Throughout this paper, $A$ always denotes an $m\times n$ matrix with real entries; we assume that $n>m$ and that ${\rm rank}(A)=m$. We denote the columns of $A$ by ${\bf a}_1,{\bf a}_2,\ldots,{\bf a}_n$. The random vector $X=(X_1,X_2,\ldots,X_n)$ is defined as in the introduction, so as to maximize the entropy $H(X)$ subject to the constraint ${\bf E}[AX]={\bf b}=(b_1,b_2,\ldots,b_m)\in {\Bbb R}^m$. We define the parameters $q_j, z_j$ as in the introduction.

We define the {\sl point concentration} of a discrete random variable $Y$ by
$${\rm conc}(Y) := \max_y {\bf Pr}[Y=y].$$
An upper bound on ${\rm conc}(AX)$ is, necessarily, also an upper bound on\break
${\bf Pr}[AX={\bf b}] = {\bf Pr}[X\in P]$. Therefore, we have
\begin{equation}
|P\cap {\Bbb Z}^n| \quad \leq \quad e^{H(X)} {\rm conc}(AX).
\end{equation}

\subsection{Results}
Under the hypotheses above, we prove:

\begin{Thm}
\begin{align*}
|P\cap {\Bbb Z}^n| \quad &\leq \quad e^{H(X)} \min_{{\bf a}_{j_1},\ldots,{\bf a}_{j_m} {\rm lin. indep.}} (1-q_{j_1})(1-q_{j_2})\cdots (1-q_{j_m}) \\ \notag
&= \quad e^{H(X)} \min_{{\bf a}_{j_1},\ldots,{\bf a}_{j_m} {\rm lin. indep.}} \prod_{i=1}^m \frac{1}{z_{j_i}+1}.
\end{align*}
\end{Thm}

\begin{Cor}
Let $I_1,I_2\ldots,I_p$ be $m$-element subsets of $\{1,2,\ldots,n\}$,
$$I_k = \{j_{k1},j_{k2},\ldots,j_{km}\},$$
such that ${\bf a}_{j_{k1}},\ldots,{\bf a}_{j_{km}}$ form a basis for ${\Bbb R}^m$ $(1\leq k\leq p)$, and such that\break
$I_1\cup I_2\cup\cdots\cup I_p = \{1,2,\ldots,n\}$. Then
$$|P\cap {\Bbb Z}^n| \quad \leq \quad e^{H(X)} \Big(\frac{1}{{\bf E}[\bar X]+1}\Big)^m,$$
where ${\bar X}$ is a geometrically distributed random variable with entropy equal to $\frac{1}{pm} H(X)$.
\end{Cor}

(A formula for the entropy of a geometrically distributed random variable is given in section 1, (1).)

\begin{Thm}
Suppose that $n=pm$ for some integer $p$, that $A$ has integer entries, and that ${\bf a}_{(k-1)m+1}, {\bf a}_{(k-1)m+2}, \ldots, {\bf a}_{km}$ are linearly independent for $1\leq k\leq p$. Assume that $\langle {\bf a}_j,{\bf b}\rangle >0$ for $1\leq j\leq n$. Define
$$q^\vee_i := \min\{q_{(k-1)m+i}: 1\leq k\leq p\} \quad\quad (1\leq i\leq m).$$

Then there exist constants $C=C(q^\vee_1,\ldots,q^\vee_m)$ and $C'=C'(q^\vee_1,\ldots,q^\vee_m)$, with $C'<1$, such that
$$|P\cap {\Bbb Z}^n| \quad \leq \quad e^{H(X)}  \big(Cp^{-m/2} + (C')^p\big).$$
\end{Thm}
(In fact, there is a one-parameter family of pairs of constants $(C,C')$ for which this statement holds. Explicit formulas and bounds for $C$ and $C'$ are provided in section 5.)

\begin{Thm}
Suppose that $n=pm$ for some integer $p$ and that, for each $i=1,2,\ldots,m$, we have
${\bf a}_i = {\bf a}_{m+i} = {\bf a}_{2m+i} = \cdots = {\bf a}_{(p-1)m+i}$, where $\{{\bf a}_1,{\bf a}_2,\ldots,{\bf a}_m\}$ is a basis for ${\Bbb R}^m$. (That is to say, the columns of $A$ cycle through a basis of ${\Bbb R}^m$ periodically.) Then
$$|P\cap {\Bbb Z}^n| \quad \lesssim \quad e^{H(X)} \prod_{i=1}^m \Big( \frac{\pi p}{6}\big(\lfloor{\bf E}(X_i)+1\rfloor^2 -1\big)\Big)^{-1/2}.$$
(Here $\lesssim$ means that, given fixed ${\bf a}_1,{\bf a}_2,\ldots,{\bf a}_m$, the expression on the left side is bounded above by a function which is asymptotic to the expression on the right side as $p\rightarrow\infty$.)
\end{Thm}

\subsection{Plan of paper}
In section 3, we discuss these results in the context of prior work, and give examples of their use. In section 4, we prove Theorem 1 and Corollary 1. The most substantial portion of the paper is section 5, in which we prove Theorem 2, then bound the constants appearing in it. In section 6, we prove Theorem 3.

\section{Discussion and examples}
The concentration of sums of random variables is a richly studied subject. The particular program of obtaining {\sl upper} bounds, sometimes called ``anti-concentration results,'' may be considered to have originated with the {\sl Littlewood--Offord problem}. This problem asked for the maximum concentration of
$$\varepsilon_1 a_1+\varepsilon_2 a_2+\cdots+\varepsilon_n a_n$$
when $a_1,a_2,\ldots,a_n$ are integers and $\varepsilon_1,\varepsilon_2,\ldots,\varepsilon_n$ are symmetric Bernoulli random variables. The exact solution, which is of order $O(n^{-1/2})$, was provided by Erd\H os~\cite{Erdos1945}.

Hal\'asz~\cite{Halasz1977} extended this result to random sums of $m$-vectors (again with symmetric Bernoulli coefficients), obtaining a bound of order $O(n^{-m/2})$ under conditions ensuring that the vectors are reasonably ``spread out'' in ${\Bbb R}^m$ (i.e., not excessively close to a proper subspace). Hal\'asz's results pertain to the {\sl small ball} concentration of $\varepsilon_1 a_1+\varepsilon_2 a_2+\cdots+\varepsilon_n a_n$, but can be specialized to point concentration. These results, which Hal\'asz proved using a Fourier-theoretic lemma of Ess\'een, were subsequently reproduced by Oskolkov~\cite[notes by Howard]{Howard} using {\sl rearrangement inequalities}. Theorem 2, herein, arrives at a similar conclusion when the Bernoulli coefficients are replaced by geometric ones. In particular, Theorem 2 implies the following Gaussian-like asymptotics:
\begin{Cor}
Suppose that a subset of the columns of $A$ can be partitioned into $p$ bases for ${\Bbb R}^m$. Then for $\min_j q_j$ bounded away from 0, the point concentration of $AX$ is $O(p^{-m/2})$ as $p\rightarrow\infty$.
\end{Cor}

Our proof of Theorem 2 hews closely to the method of~\cite{Howard}. For other approaches to anti-concentration inequalities, see~\cite{RV2010},~\cite{TV2009}.

Theorem 2 is essentially an asymptotic result; although we give explicit formulas for $C$ and $C'$, the bounds obtained from Theorem 2 are typically only strong when $p$ is large, i.e., when $n\gg m$. (For further remarks on this theme, see the end of section 5.1.) By contrast, Theorem 1 and its corollary are non-asymptotic, and are apparently most effective when $n \not\gg m$. They are also relatively straightforward, but do not capture the $O(p^{-m/2})$ behavior of ${\rm conc}(AX)$. Thus, Theorem 1 and Theorem 2 may be seen as filling somewhat different niches. Theorem 3 gives a more ideal bound, combining all the attractive features of Theorems 1 and 2, but is pertinent only to a very special case (the easiest, where a local central limit theorem is available). Relying for its proof on notions from the theory of partially ordered sets, Theorem 3 may serve as a suggestion of how combinatorics can be brought to bear on this problem.

\subsection{Examples}
Given nonnegative vectors $R\in{\Bbb R}^r,~~C\in{\Bbb R}^s$, the {\sl transportation polytope} $\Pi(R,C)$ is defined as the set of all nonnegative $r\times s$ matrices whose row sums and column sums are the coordinates of $R$ and $C$, respectively. Such a matrix with integer entries is called a {\sl contingency table}.

We may use Theorem 1 to bound the number of $4\times 4$ contingency tables with given ``margins'' $R$ and $C$. For example, let $R=(108,286,71,127)$ and $C=(220,215,93,64)$, as in a table appearing in~\cite{DE1985} which has become a standard benchmark in the literature on contingency tables. The actual number of tables with these row and column sums is $1.23\times 10^{15}$. (It can be computed exactly, as the dimension is fairly low: the defining matrix $A$ for $\Pi(R,C)$ is $7\times 16$.)

Let $X$ be a random matrix taking the maximum-entropy distribution on ${\Bbb Z}_{\geq 0}^{4\times 4}$, under the constraint that ${\bf E}[X]\in\Pi(R,C)$. Solving the convex optimization problem described in section 1, we compute
\[
{\bf E}[X] =
\left( {\begin{matrix}
 36.4 & 36.0 & 20.6 & 14.9  \\
 117.2 & 113.4 & 34.3 & 21.2  \\
 22.2 & 22.0 & 15.1 & 11.7 \\
 44.2 & 43.6 & 23.0 & 16.2 \\
 \end{matrix} } \right)
\]
and $H(X)=2.96\times 10^{30}$. Theorem 1 then gives
\begin{align*}
|P\cap {\Bbb Z}^n| &\leq \frac{2.96\times 10^{30}}{(1+36.4)(1+117.2)(1+113.4)(1+34.3)(1+21.2)(1+22.2)(1+44.2)} \\
&= 7.14\times 10^{18},
\end{align*}
off by a factor of about 5800. Computation of similar examples suggests that the relative error depends mainly on the dimensions of $R$ and $C$, and not on the magnitude of their entries.

Theorem 2 performs relatively poorly in these examples, but is much more effective than Theorem 1 when $n$ is large compared to $m$. For instance, consider the simplex
$$\Sigma^n(r) := \big\{(x_1,\ldots,x_n):\quad x_1,\ldots,x_n\geq 0,\quad \|{\bf x}\|_1 = r\big\},$$
which has ${n+r-1\choose r}$ integer points. Let $0<\delta<\frac{1}{2}$. Then, choosing $\gamma=\frac{\pi r^\delta}{\sqrt{n}}$ in the statement of Theorem 2a (see section 5), one obtains as a conclusion an upper bound on $|\Sigma^n(r)\cap{\Bbb Z}^n|$ which is precisely asymptotic to ${n+r-1\choose r}$ as $n\rightarrow\infty$, if $r$ grows as $\Theta(n^\varepsilon)$ for some $\varepsilon\in(0,1)$.

For $r=10$ and $n=1000$, the optimal result of Theorem 2a (achieved when $\gamma=0.172$) is  an upper bound of $3.14\times 10^{23}$, which may be compared with an exact count of $2.88\times 10^{23}$ integer points. By comparison, when $r=100$ and $n=10000$, the optimal result of Theorem 2a (achieved when $\gamma=0.0645$) is an upper bound of $1.774\times 10^{242}$ integer points; the exact count is $1.755\times 10^{242}$, and the relative error is about $1.1\%$.

\section{Proofs of Theorem 1 and Corollary 1}

We prove Theorem 1 by means of the following simple fact:
\begin{Lem} If $X, Y$ are independent, discrete random variables, then ${\rm conc}(X+Y)\leq {\rm conc}(X)$. \end{Lem}
{\sl Proof.}\quad Observe that ${\rm conc}(X+Y)$ is a weighted average of values of the probability mass function of $X$, of which the largest is ${\rm conc}(X)$. $\square$

{\sl Proof of Theorem 1.}\quad Using Lemma 1 and the previously mentioned properties of geometric random variables,
\begin{align*}
{\rm conc}(X_1{\bf a}_1+\cdots+X_n{\bf a}_n) \quad &\leq \quad \min_{{\bf a}_{j_1},\ldots,{\bf a}_{j_m} {\rm lin. indep.}} {\rm conc}(X_{j_1}{\bf a}_{j_1}+\cdots+X_{j_m}{\bf a}_{j_m}) \\
&\leq \quad \min_{{\bf a}_{j_1},\ldots,{\bf a}_{j_m} {\rm lin. indep.}} {\bf Pr}[X_{j_1}=\cdots =X_{j_m}=0] \\
&= \quad \min_{{\bf a}_{j_1},\ldots,{\bf a}_{j_m} {\rm lin. indep.}} (1-q_{j_1})(1-q_{j_2})\cdots (1-q_{j_m}) \\
&= \quad \min_{{\bf a}_{j_1},\ldots,{\bf a}_{j_m} {\rm lin. indep.}} \prod_{i=1}^m \frac{1}{z_{j_i}+1}.
\end{align*}
By section 2.1, (2), it follows that
\begin{align*}
|P\cap {\Bbb Z}^n| \quad &\leq \quad e^{H(X)} \min_{{\bf a}_{j_1},\ldots,{\bf a}_{j_m} {\rm lin. indep.}} (1-q_{j_1})(1-q_{j_2})\cdots (1-q_{j_m}) \\ \notag
&= \quad e^{H(X)} \min_{{\bf a}_{j_1},\ldots,{\bf a}_{j_m} {\rm lin. indep.}} \prod_{i=1}^m \frac{1}{z_{j_i}+1}. \quad\quad\blacksquare
\end{align*}

To prove Corollary 1, we will require this fact whose proof is deferred until after the proof of Corollary 1:
\begin{Lem}
Among all vectors $Y := (Y_1,Y_2,\ldots,Y_m)$ of independent, geometrically distributed random variables with fixed joint entropy $\Omega$, the highest concentration ${\rm conc}(Y)$ is achieved when $Y_1,Y_2,\ldots,Y_m$ are identically distributed.
\end{Lem}

{\sl Proof of Corollary 1.}\quad For $I\subset \{1,2,\ldots,n\}$, let $H(X_I)$ denote the joint entropy of $\{X_j: j\in I\}$. Since $X_1,\ldots,X_n$ are pairwise independent, we have $H(X_I) = \sum_{j\in I} H(X_j)$.

Since the sets $I_1,I_2,\ldots,I_p$ cover $\{1,2,\ldots,n\}$, we have
$$H(X)\leq \sum_{k=1}^p H(X_{I_k}),$$
and thus by the pigeonhole principle
$$H(X_{I_k}) \geq \frac{1}{p} H(X)$$
for some $k\in\{1,\ldots,p\}$.
By Lemma 2, the concentration of the vector $(X_{j_{k1}},\ldots,X_{j_{km}})$ is maximized when $X_{j_{k1}},\ldots,X_{j_{km}}$ are identically distributed. In this case, each has entropy equal to $\frac{1}{m} H(X_{I_k})$, which is greater than or equal to $H(\bar X) = \frac{1}{pm} H(X)$; we pause to note that the entropy and the expectation of a geometric random variable are monotonically increasing functions of one another. Thus (as in the proof of Theorem 1),
\begin{align*}
{\rm conc}(AX) \quad&\leq\quad {\rm conc}(X_{j_{k1}}{\bf a}_{j_{k1}} + \cdots + X_{j_{km}}{\bf a}_{j_{km}}) \\
&\leq\quad \Big(\frac{1}{{\bf E}[\bar X]+1}\Big)^m,
\end{align*}
so Corollary 1 follows by section 2.1, (2). $\blacksquare$

{\sl Proof of Lemma 2.}\quad
Since $Y_i$ is geometrically distributed $(1\leq i\leq m)$, there exist parameters $r_i\in[0,1)$ such that
$${\bf Pr}[Y_i=k] = (1-r_i)r_i^k \quad{\rm for}\quad k\in{\Bbb Z}_{\geq 0}.$$
The concentration of $Y$ is $\prod_{i=1}^m (1-r_i)$, so we must show that this expression is maximized (for fixed $\Omega$) when $r_1 = \ldots = r_m$.

We introduce the changes of variable $s_i:=\frac{1}{1-r_i}$,\quad $t_i:=\ln s_i$. (Thus $1-r_i=\frac{1}{s_i}$, and $s_i=e^{t_i}$, where $t_i\in[0,\infty)$.) Also, let
$$\omega(t) := (1-e^t)\ln(1-e^{-t})+t.$$
Now
\begin{align*}
\Omega &= \sum_{i=1}^m \frac{r_i}{1-r_i}\ln\frac{1}{r_i} + \ln\frac{1}{1-r_i} \\
&= \sum_{i=1}^m (s_i-1) \ln\frac{s_i}{s_i-1} + \ln s_i \\
&= \sum_{i=1}^m (e^{t_i}-1) \ln\frac{e^{t_i}}{e^{t_i}-1} + t_i \\
&= \sum_{i=1}^m (1-e^{t_i}) \ln (1-e^{-t_i}) + t_i \\
&= \sum_{i=1}^m \omega(t_i),
\end{align*}
and
$$\prod_{i=1}^m (1-r_i) = \exp\bigg(-\sum_{i=1}^m t_i\bigg).$$

The following three statements are equivalent:
\begin{enumerate}
\item For $\Omega$ fixed, $\prod\limits_i (1-r_i)$ is maximized when $r_1=\cdots=r_m$.
\item For $\Omega$ fixed, $\sum\limits_i t_i$ is minimized when $t_1=\ldots=t_m$.
\item If $\sum\limits_i t_i$ is fixed and $\Omega$ free to vary, then $\Omega$ is maximized when $t_1=\ldots=t_m$.
\end{enumerate}

The equivalence of statements (1) and (2) is clear. To see that (2) and (3) are equivalent, it is enough to observe that $\Omega$ is increasing with respect to each of $t_1,\ldots,t_m$. Thus to prove (1), which is the assertion of the lemma, it will suffice for us to prove (3).

Writing $s:=e^t$, we obtain
\begin{align*}
{d\omega\over dt} &= (1-e^t)\bigg({e^{-t}\over 1-e^{-t}}\bigg) - e^t\ln(1-e^{-t}) + 1 \\
&= -e^t\ln(1-e^{-t})
\end{align*}
and
\begin{align*}
{d^2\omega\over dt^2} &= -e^t\cdot {e^{-t}\over 1-e^{-t}} - e^t\ln(1-e^{-t}) \\
&= -{1\over 1-{1\over s}} - s\ln\bigg(1-{1\over s}\bigg) \\
&= -{s\over s-1} + s\ln {s\over s-1} \\
&= -s\bigg({1\over s-1}\bigg) + s\ln\bigg(1+{1\over s-1}\bigg) \\
&\leq 0,
\end{align*}
since $\ln(1+x)\leq x$ for $x\geq 0$. This shows that $\omega(t)$ is concave for $t\geq 0$, which implies (3) and so completes the proof of the lemma. $\blacksquare$

\section{Proof of Theorem 2}

We begin by restating the theorem with explicit formulas for all constants:
\begin{Thm2a}
Assume the definitions and notation from section 2.1.

Suppose that $n=pm$ for some integer $p$, that $A$ has integer entries, and that ${\bf a}_{(k-1)m+1}, {\bf a}_{(k-1)m+2}, \ldots, {\bf a}_{km}$ are linearly independent for $1\leq k\leq p$. Assume that $\langle {\bf a}_j,{\bf b}\rangle >0$ for $1\leq j\leq n$. Let $\gamma>0$. Define constants
\begin{align*}
\alpha_j:=\frac{2q_j}{(1-q_j)^2} &\quad\quad (1\leq j\leq n), \\
\alpha^\vee_i := \min\{\alpha_{(k-1)m+i}: 1\leq k\leq p\} &\quad\quad (1\leq i\leq m), \\
q^\vee_i := \min\{q_{(k-1)m+i}: 1\leq k\leq p\} &\quad\quad (1\leq i\leq m), \\
c_i := \max\bigg\{\frac{1}{\gamma^2}\ln\Big[1+\alpha^\vee_i \big(1-\cos\frac{\gamma}{\sqrt{\alpha^\vee_i}}\big)\Big],\quad
\frac{1}{\alpha^\vee_i\pi^2}\ln\Big[1+2\alpha^\vee_i\Big] \bigg\} &\quad\quad (1\leq i\leq m), \\
C := \prod\limits_{i=1}^m (2\pi c_i\alpha^\vee_i)^{-1/2}, \\
C' := \max\limits_{1\leq i\leq m} e^{-\gamma^2 c_i/2}.
\end{align*}
Then
$$|P\cap {\Bbb Z}^n| \quad \leq \quad e^{H(X)}  \big(Cp^{-m/2} + (C')^p\big).$$
\end{Thm2a}

All notation introduced in Theorem 2a is used throughout this section, and all its hypotheses (importantly, the integrality of $A$) are assumed to hold. In subsection 5.1, we introduce a series of definitions and lemmas, then prove Theorem 2a under assumption of the lemmas. In subsection 5.2, we prove the lemmas in turn. For bounds on the constants $C$ and $C'$, see subsection 5.3.

\subsection{Supporting results and proof of Theorem 2a}
\begin{Def}
For $1\leq k\leq p$, define the function $\Pi_k: {\Bbb R}^m\rightarrow{\Bbb R}$ by
\begin{align*}
\Pi_k({\bf t}) := \prod_{j=(k-1)m+1}^{km} \frac{1} {\sqrt{1+\alpha_j(1-\cos\langle{\bf t}, {\bf a}_j\rangle)}} \quad\quad &{\rm for}~{\bf t}\in (-\pi,\pi]^m, \\
\Pi_k({\bf t}) := 0 \quad\quad &{\rm for}~{\bf t}\not\in(-\pi,\pi]^m.
\end{align*}
\end{Def}

\begin{Lem}
Given the definition above,
$${\rm conc}(AX) \leq \frac{1}{(2\pi)^m} \int_{(-\pi,\pi]^m} \Pi_1 \Pi_2 \cdots \Pi_p ~d{\bf t}.$$
\end{Lem}

\begin{Def}
Given a measurable function $\Phi:{\Bbb R}^m\rightarrow{\Bbb R_{\geq 0}}$, we define its {\rm epigraphs}
$$\Gamma_{\geq\tau} (\Phi) := \{{\bf t}\in{\Bbb R}^m : \Phi({\bf t}) \geq\tau \}$$
for all $\tau>0$.

Suppose $\Phi$ {\rm vanishes at infinity}, meaning that $\Gamma_{\geq\tau} (\Phi)$ has finite volume for each $\tau>0$. Then we define its {\rm symmetrically decreasing rearrangement} as the function $\Phi^*:{\Bbb R}^m\rightarrow{\Bbb R_{\geq 0}}$ given by
$$\Phi^*({\bf t}) := \max\Big\{\tau:~~{\rm vol}\big(\Gamma_{\geq\tau}(\Phi)\big)\geq \|{\bf t}\|^m v_m\Big\},$$
where $v_m$ denotes the volume of the unit ball in ${\Bbb R}^m$.
\end{Def}

The theory of symmetrically decreasing rearrangements is treated in~\cite{Burchard}, and we do not develop it fully here. The important properties of $\Phi^*$ are that
\begin{itemize}
\item $\Phi^*$ is symmetrically decreasing, i.e., $\|{\bf t}\| \geq \|{\bf s}\| \Rightarrow \Phi^*({\bf t}) \leq \Phi^*({\bf s})$; and
\item $\Phi^*$ is equimeasurable with $\Phi$, i.e., ${\rm vol}(\Gamma_{\geq\tau} (\Phi^*)) = {\rm vol}(\Gamma_{\geq\tau} (\Phi))$ for all $\tau>0$.
\end{itemize}
Note that $\Phi^*$ is the unique function with these properties, up to a difference on a set of measure zero.

\begin{Lem}
Given the definition above,
$$\int_{(-\pi,\pi]^m} \Pi_1 \Pi_2 \cdots \Pi_p ~d{\bf t} \leq \int_{{\Bbb R}^m} \Pi^*_1 \Pi^*_2 \cdots \Pi^*_p ~d{\bf t}.$$
\end{Lem}

\begin{Def}
For $1\leq k\leq p$, define the function $\Pi^{\rm rect}_k: {\Bbb R}^m\rightarrow{\Bbb R}$ by
\begin{align*}
\Pi^{\rm rect}_k({\bf t}) := \prod_{i=1}^{m} \frac{1} {\sqrt{1+\alpha_{(k-1)m+i}(1-\cos{\bf t}_i)}}
\quad\quad &{\rm for}~{\bf t}\in (-\pi,\pi]^m, \\
\Pi^{\rm rect}_k({\bf t}) := 0 \quad\quad &{\rm for}~{\bf t}\not\in(-\pi,\pi]^m.
\end{align*}
\end{Def}

The formula for $\Pi^{\rm rect}_k$ differs from that for $\Pi_k$ in that the linear form $\langle{\bf t}, {\bf a}_{(k-1)m+i}\rangle$ in the denominator of $\Pi_k$ is replaced by ${\bf t}_i$. Effectively, each basis
$${\bf a}_{(k-1)m+1}, {\bf a}_{(k-1)m+2}, \ldots, {\bf a}_{km}$$
of ${\Bbb R}^m$ is replaced by a standard basis. This will make $\Pi^{\rm rect}_k$ easier to work with than $\Pi_k$.

\begin{Lem}
Let $1\leq k\leq p$. Then
$${\rm vol}\big(\Gamma_{\geq\tau} (\Pi^{\rm rect}_k)\big) = {\rm vol}\big(\Gamma_{\geq\tau} (\Pi_k)\big)$$
for all $\tau>0$, and $(\Pi^{\rm rect}_k)^* \equiv \Pi^*_k$.
\end{Lem}

\begin{Lem}[Isotonicity of rearrangement]
Suppose $\Phi,\Psi:{\Bbb R}^m\rightarrow{\Bbb R_{\geq 0}}$ are measurable functions vanishing at infinity. Let $\tau$ denote a constant. Then:
\begin{enumerate}
\item If $\Phi({\bf t}) \geq \Psi({\bf t})$ for all ${\bf t}$, then $\Phi^*({\bf t}) \geq \Psi^*({\bf t})$ for all ${\bf t}$.
\item If $\Phi({\bf t}) \geq \max\{\Psi({\bf t}), \tau\}$ for all ${\bf t}$, then $\Phi^*({\bf t}) \geq \max\{\Psi^*({\bf t}), \tau\}$ for all ${\bf t}$.
\end{enumerate}
\end{Lem}

\begin{Lem}
Define $\alpha^\vee_i$ and $c_i$ as in the statement of Theorem 2a.

Then, for $0\leq t\leq \min\Big\{\frac{\gamma}{\sqrt{\alpha^\vee_i}},~~\pi\Big\}$, we have
$1+\alpha^\vee_i(1-\cos t) \geq e^{c_i\alpha^\vee_i t^2}$.
\end{Lem}

\begin{Lem}
For each $k=1,2,\ldots,p$, and for all ${\bf t}\in{\Bbb R}^m$, we have
$$\Pi^{\rm rect}_k({\bf t}) \leq \max\bigg\{\prod_{i=1}^m e^{-c_i\alpha^\vee_i{\bf t}_i^2/2},\quad C'\bigg\}.$$
\end{Lem}

Given the above lemmas, we can prove Theorem 2a:

{\sl Proof of Theorem 2a.}\quad Using Lemmas 3, 4, and 5, we have
\begin{align*}
{\rm conc}(AX) &\leq \frac{1}{(2\pi)^m} \int_{(-\pi,\pi]^m} \Pi_1 \Pi_2 \cdots \Pi_p ~d{\bf t} \\
&\leq \frac{1}{(2\pi)^m} \int_{{\Bbb R}^m} \Pi^*_1 \Pi^*_2 \cdots \Pi^*_p ~d{\bf t} \\
&= \frac{1}{(2\pi)^m} \int_{{\Bbb R}^m} (\Pi^{\rm rect}_1)^* (\Pi^{\rm rect}_2)^* \cdots (\Pi^{\rm rect}_p)^* ~d{\bf t}.
\end{align*}
We may instead take either of the last two integrals over $B$, the closed ball of volume $(2\pi)^m$ centered at the origin in ${\Bbb R}^m$, since the integrands are zero outside this ball.

By Lemmas 6 and 8, we have
\begin{align*}
\frac{1}{(2\pi)^m} \int_B ~~& (\Pi^{\rm rect}_1)^* (\Pi^{\rm rect}_2)^* \cdots (\Pi^{\rm rect}_p)^* ~d{\bf t}
\quad \\
&\leq \quad
\frac{1}{(2\pi)^m} \int_B ~~\prod_{k=1}^p \Bigg( \max\bigg\{ \Big( \prod_{i=1}^m e^{-c_i \alpha^\vee_i {\bf t}_i^2 /2} \Big)^*, \quad C' \bigg\} \Bigg) ~d{\bf t} \\
\quad &= \quad \frac{1}{(2\pi)^m} \int_B ~~\Bigg( \max\bigg\{ \Big( \prod_{i=1}^m e^{-c_i \alpha^\vee_i {\bf t}_i^2 /2} \Big)^*, \quad C' \bigg\} \Bigg)^p ~d{\bf t} \\
\quad &= \quad \frac{1}{(2\pi)^m} \int_{(-\pi,\pi]} \Bigg( \max\bigg\{ \prod_{i=1}^m e^{-c_i \alpha^\vee_i {\bf t}_i^2 /2}, \quad C' \bigg\} \Bigg)^p ~d{\bf t}.
\end{align*}

This last integral is bounded above by
\begin{align*}
&\quad \frac{1}{(2\pi)^m} \Bigg[ \int_{{\Bbb R}^m}~ \bigg( \prod_{i=1}^m e^{-c_i \alpha^\vee_i {\bf t}_i^2 /2} \bigg)^p ~d{\bf t} \quad+\quad \int_B~ (C')^p~d{\bf t} \Bigg] \\
= &\quad \frac{1}{(2\pi)^m} \Bigg[ \int_{{\Bbb R}^m}~ \exp\bigg( -p \sum_{i=1}^m -c_i \alpha^\vee_i {\bf t}_i^2 /2 \bigg) ~d{\bf t} \quad+\quad (2\pi)^m (C')^p~d{\bf t} \Bigg] \\
= &\quad \frac{1}{(2\pi)^m} \cdot (2\pi)^{m/2} p^{-m/2} \prod_{i=1}^m (c_i\alpha^\vee_i)^{-1/2} \quad+\quad (C')^p \\
= &\quad Cp^{-m/2} + (C')^p.
\end{align*}

Now, a technical remark. In integrating the Gaussian term, we assumed $c_i\alpha^\vee_i>0$. To see why this is necessarily true, note that we assumed, in the statement of Theorem 2a, that $\langle {\bf a}_j,{\bf b}\rangle >0$ for $1\leq j\leq n$. Thus $P$ is not contained in any coordinate hyperplane of ${\Bbb R}^n$. Recall section 1, (1), which gives the entropy $H(X)$ in terms of $z_1,z_2,\ldots,z_n$ (the coordinates of ${\bf E}[X])$. One may check that $\frac{\partial}{\partial z_j} H(X)=\infty$ when $z_j=0$, but is finite when $z_j>0$. Therefore the maximum-entropy distribution for $X$ does not take expected value on a coordinate hyperplane; therefore,\break
$c_i\alpha^\vee_i>0$.

Theorem 2a now follows by section 2.1, (2). $\blacksquare$

{\sl Remarks.} Our strategy for bounding ${\rm conc}(AX)$, carried out above, may be motivated as follows. First, we obtain an integral formula for the probability mass function of $AX$, derived from its Fourier transform (Lemma 3). The integrand splits into $n$ factors, which we then group into maximal subproducts such that the factors in each subproduct behave like independent random variables on the domain of integration. The worst case is now that these subproducts themselves are ``completely non-independent,'' that is, that they decay identically; this is the significance of Lemmas 4 and 5, and of the definitions of $q^\vee_i$ and $\alpha^\vee_i$. We bound the decay of the integrand near the origin by a Gaussian (Lemma 8), explaining the appearance of the $Cp^{-m/2}$ term in the conclusion of Theorem 2a. Away from the origin, we simply bound each subproduct by the constant $C'$, giving the $(C')^p$ term. The parameter $\gamma$ controls the boundary between the two approximation regimes.

This two-regime bound (with arbitrary parameter $\gamma$) is sufficient for Corollary 2, as the $(C')^p$ term is asymptotically negligible as $p\rightarrow\infty$. However, for non-asymptotic computations, the crudity of the approximation away from the origin is very noticeable. The $(C')^p$ term can be replaced by a more sensitive approximation, at the cost of simplicity. We do not pursue this goal here.

\subsection{Proofs of preceding lemmas}~

{\sl Proof of Lemma 3.}\quad 
In~\cite{BH2010}, Lemma 8.1, the following integral representation is proved:
$${\bf Pr}[AX={\bf b}] \quad=\quad \frac{1}{(2\pi)^m} \int_{(-\pi,\pi]^m} e^{-i\langle {\bf t},{\bf b}\rangle} \prod_{j=1}^n \frac{1-q_j} {1-q_j e^{i\langle{\bf t}, {\bf a}_j\rangle}} ~d{\bf t},$$
where ${\bf b}$ is an arbitrary ${\Bbb Z}_{\geq 0}$-vector.
It follows that
\begin{align*}
{\rm conc}(AX) \quad&\leq\quad \frac{1}{(2\pi)^m} \int_{(-\pi,\pi]^m} \Bigg |e^{-i\langle {\bf t},{\bf b}\rangle} \prod_{j=1}^n \frac{1-q_j} {1-q_j e^{i\langle{\bf t}, {\bf a}_j\rangle}} \Bigg | ~d{\bf t} \\
\quad&=\quad \frac{1}{(2\pi)^m} \int_{(-\pi,\pi]^m} \prod_{j=1}^n \frac{1-q_j} {\sqrt{1+q_j^2-2q_j\cos\langle{\bf t}, {\bf a}_j\rangle}} ~d{\bf t} \\
\quad&=\quad \frac{1}{(2\pi)^m} \int_{(-\pi,\pi]^m} \Pi_1 \Pi_2 \cdots \Pi_p ~d{\bf t},
\end{align*}
where the last two steps are straightforward simplification. $\square$

{\sl Proof of Lemma 4.}\quad
The {\sl Hardy--Littlewood inequality}~\cite{Burchard} states that for measurable functions $\Phi,\Psi:{\Bbb R}^m\rightarrow{\Bbb R_{\geq 0}}$ vanishing at infinity, one has
$$\int_{{\Bbb R}^m} \Phi({\bf t}) \Psi({\bf t}) ~d{\bf t} \leq \int_{{\Bbb R}^m} \Phi^*({\bf t}) \Psi^*({\bf t}) ~d{\bf t},$$
provided that the integral on the right-hand side converges. Thus we obtain
$$\int_{(-\pi,\pi]^m} \Pi_1 \Pi_2 \cdots \Pi_p ~d{\bf t} \leq \int_{{\Bbb R}^m} \Pi^*_1 \Pi^*_2 \cdots \Pi^*_p ~d{\bf t}$$
by induction on $p$. $\square$

{\sl Proof of Lemma 5.}

Let $A_*$ be the $m\times m$ matrix whose rows are ${\bf a}_{(k-1)m+1}^T, {\bf a}_{(k-1)m+2}^T, \ldots, {\bf a}_{km}^T$, and define ${\mathcal A_*}:{\Bbb R}^m\rightarrow{\Bbb R}^m$ as the linear map ${\bf t}\mapsto A_*{\bf t}$. Thus,
$${\mathcal A_*}({\bf t})_i = \langle{\bf t},{\bf a}_{(k-1)m+i}\rangle\quad(1\leq i\leq m).$$
This map ${\mathcal A_*}$ scales the volume of measurable sets uniformly by a factor of $d:=|\det(A_*)|$, and takes the lattice $\Lambda := (2\pi{\Bbb Z})^m$ to the lattice
$$\Lambda' := 2\pi{\Bbb Z} [{\bf col}_1(A_*), {\bf col}_2(A_*), \ldots, {\bf col}_m(A_*)].$$

Let $K := (-\pi,\pi]^m$ and let $K' := {\mathcal A_*}(K)$. Since $K$ is a fundamental region of $\Lambda$, it follows that $K'$ is a fundamental region of $\Lambda'$. Moreover, we assumed $A$ to have integer entries, so $\Lambda'$ is a sublattice of index $d$ in $\Lambda$, and the induced map of tori $\phi: {\Bbb R}^m/\Lambda' \rightarrow {\Bbb R}^m/\Lambda$ is an even covering of order $d$.

Identifying $K$ with ${\Bbb R}^m/\Lambda$ and $K'$ with ${\Bbb R}^m/\Lambda'$, we may regard $\phi$ as a map from $K'$ to $K$, and $\phi\circ{\mathcal A_*}$ as a self-map of $K$. If $U\subseteq K$ is a measurable set, then $(\phi\circ{\mathcal A_*})^{-1}(U)$ is the union of $d$ disjoint preimages each of volume $\frac{{\rm vol}(U)}{d}$. Thus, ${\rm vol}((\phi\circ{\mathcal A_*})^{-1}(U)) = {\rm vol}(U)$.

Observe that $\cos{\bf t}_i = \cos(\phi({\bf t})_i)$ for all ${\bf t}$. Therefore
\begin{align*}
\Gamma_{\geq\tau} (\Pi_k) &= {\mathcal A_*}^{-1} (\Gamma_{\geq\tau} (\Pi^{\rm rect}_k)) \\
&= (\phi\circ{\mathcal A_*})^{-1} (\Gamma_{\geq\tau} (\Pi^{\rm rect}_k))
\end{align*}
from which it follows that
$${\rm vol}\big(\Gamma_{\geq\tau} (\Pi^{\rm rect}_k)\big) = {\rm vol}\big(\Gamma_{\geq\tau} (\Pi_k)\big).$$
This conclusion holds for all $\tau>0$, so it follows from the definition of the symmetrically decreasing rearrangement that $(\Pi^{\rm rect}_k)^* \equiv \Pi^*_k$. $\square$

{\sl Proof of Lemma 6.}\quad
We prove (1) by contradiction. Suppose that $\Phi({\bf t}) \geq \Psi({\bf t})$ for all ${\bf t}$, but suppose $\Phi^*({\bf t}_0) < \Psi^*({\bf t}_0)$ for some ${\bf t}_0$. Let $\tau_0:= \Psi^*({\bf t}_0)$. Then
$${\rm vol}\big(\Gamma_{\geq\tau_0}(\Phi)\big) < \|{\bf t}_0\|^m v_m \leq {\rm vol}\big(\Gamma_{\geq\tau_0}(\Psi)\big),$$
where $v_m$ is the volume of the unit ball in ${\Bbb R}^m$. It follows that $\Gamma_{\geq\tau_0}(\Psi) \backslash \Gamma_{\geq\tau_0}(\Phi)$ has positive measure, contradicting our assumption that $\Phi({\bf t}) \geq \Psi({\bf t})$ for all ${\bf t}$.

Statement (2) follows from (1) by the observation that $\max\{\Psi^*({\bf t}), \tau\}$ {\sl is} the symmetrically decreasing rearrangement of $\max\{\Psi({\bf t}), \tau\}$. $\square$

{\sl Proof of Lemma 7.}\quad
Recall that
$$c_i := \max\bigg\{\frac{1}{\gamma^2}\ln\Big[1+\alpha^\vee_i \big(1-\cos\frac{\gamma}{\sqrt{\alpha^\vee_i}}\big)\Big],\quad
\frac{1}{\alpha^\vee_i\pi^2}\ln\Big[1+2\alpha^\vee_i\Big] \bigg\}.$$
In particular,
$$c_i = \frac{1}{\gamma^2}\ln\Big[1+\alpha^\vee_i \big(1-\cos\frac{\gamma}{\sqrt{\alpha^\vee_i}}\big)\Big] \quad\quad {\rm if}~\alpha^\vee_i\geq \frac{\gamma^2}{\pi^2},$$
and
$$c_i = \frac{1}{\alpha^\vee_i\pi^2}\ln\Big[1+2\alpha^\vee_i\Big] \quad\quad {\rm if}~\alpha^\vee_i\leq \frac{\gamma^2}{\pi^2}.$$

Define $t_0 := \min\Big\{\frac{\gamma}{\sqrt{\alpha^\vee_i}},~~\pi\Big\}$, and define $f(t) := 1+\alpha^\vee_i(1-\cos t) - e^{c_i\alpha^\vee_i t^2}$ for $-t_0\leq t\leq t_0$.

Note that $f(0)=0$. Also, we claim that $f(t_0)=0$. This must be verified in two cases, according to whether $\alpha^\vee_i\geq \frac{\gamma^2}{\pi^2}$ or $\alpha^\vee_i\leq \frac{\gamma^2}{\pi^2}$.

If $\alpha^\vee_i\geq \frac{\gamma^2}{\pi^2}$, then $t_0=\frac{\gamma}{\sqrt{\alpha^\vee_i}}$, so
\begin{align*}
f(t_0) &= 1+\alpha^\vee_i\Big(1-\cos \frac{\gamma}{\sqrt{\alpha^\vee_i}}\Big) - \exp\Big(\frac{\alpha^\vee_i}{\gamma^2} \cdot \ln\Big[1+\alpha^\vee_i \big(1-\cos\frac{\gamma}{\sqrt{\alpha^\vee_i}}\big)\Big] \cdot \frac{\gamma^2}{\alpha^\vee_i}\Big) \\
&= 0.
\end{align*}
If $\alpha^\vee_i\leq \frac{\gamma^2}{\pi^2}$, then $t_0=\pi$, and
$$f(t_0) \quad=\quad 1+2\alpha^\vee_i - \exp\Big(\frac{1}{\alpha^\vee_i\pi^2} \cdot \ln\Big[1+2\alpha^\vee_i\Big] \cdot \alpha^\vee_i\pi^2\Big)
\quad=\quad 0.$$
This proves the claim that $f(t_0)=0$. It follows that the average value of $f'(t)$ on $[0,t_0]$ is zero.

Finally, we observe that $f'(0)=0$, and that $f(t)$ has nonpositive third derivative on $[0,t_0]$ (indeed, on $[0,\pi]$). The verification of these claims is routine and is omitted. We infer that either $f'(t)\equiv 0$ on $[0,t_0]$, or $f''(t)$ has exactly one sign change on $[0,t_0]$, from positive to negative. In the latter case, $f'(t)$ must also have exactly one sign change on $[0,t_0]$ (also from positive to negative), since its average value on the interval is zero. It follows in either case that $f(t)\geq 0$ on $[0,t_0]$, and thus on $[-t_0,t_0]$ (since $f(t)$ is an even function). This proves the lemma. $\square$

Lemma 7 is used to establish Lemma 8.

{\sl Proof of Lemma 8.}\quad 
Let
$$K:=\bigg\{{\bf t}\in{\Bbb R}^m: |{\bf t}_i|\leq \min\Big\{\frac{\gamma}{\sqrt{\alpha^\vee_i}},~~\pi\Big\}~{\rm for}~i=1,2,\ldots,m\bigg\}.$$
If ${\bf t}\in K$, then by Lemma 7,
\begin{align*}
\Pi^{\rm rect}_k({\bf t}) &= \prod_{i=1}^{m} \frac{1} {\sqrt{1+\alpha_{(k-1)m+i}(1-\cos{\bf t}_i)}} \\
&\leq \prod_{i=1}^{m} \frac{1} {\sqrt{1+\alpha^\vee_i (1-\cos{\bf t}_i)}} \\
&\leq \prod_{i=1}^{m} e^{-c_i\alpha^\vee_i{\bf t}_i^2/2}.
\end{align*}
Now suppose ${\bf t}\not\in K$. Thus, there exists some $i$ such that ${\bf t}_i > \min\Big\{\frac{\gamma}{\sqrt{\alpha^\vee_i}},~~\pi\Big\}$.\break
If ${\bf t}_i >\pi$, then we trivially have $\Pi^{\rm rect}_k({\bf t})=0\leq C'$.

Otherwise, we have ${\bf t}_i > \frac{\gamma}{\sqrt{\alpha^\vee_i}}$, and therefore
\begin{align*}
\Pi^{\rm rect}_k({\bf t}) &\leq \frac{1} {\sqrt{1+\alpha^\vee_i (1-\cos{\bf t}_i)}} \\
&\leq \frac{1} {\sqrt{1+\alpha^\vee_i \big(1-\cos \big(\gamma / \sqrt{\alpha^\vee_i}\big)\big)}} \\
&= e^{-\gamma^2 c_i/2} \\
&\leq C'.
\end{align*}

Thus whether ${\bf t}\in K$ or ${\bf t}\not\in K$, we have
$$\Pi^{\rm rect}_k({\bf t}) \leq \max\bigg\{\prod_{i=1}^m e^{-c_i\alpha^\vee_i{\bf t}_i^2/2},\quad C'\bigg\},$$
proving the lemma. $\square$

\subsection{Upper bounds on $C$, $C'$}
We now obtain
\begin{Thm2b} Defining all constants as in the statement of Theorem 2a,
$$C\leq \Bigg[\frac{\gamma}{2\sqrt{\pi\ln\big(1+\frac{2\gamma^2}{\pi^2}\big)}} \Bigg]^m \prod_{i=1}^m \frac{1-q^\vee_i}{\sqrt{q^\vee_i}}$$
and
$$C' \leq \frac{1}{\sqrt{1+\frac{2\gamma^2}{\pi^2}}}.$$
\end{Thm2b}

{\sl Remarks.}\quad Notice that as $\gamma\rightarrow\infty$, all other inputs being fixed, we have $C=O\Big(\big(\frac{\gamma}{\ln\gamma}\big)^m\Big)$ and $C'=O\big(\frac{1}{\gamma}\big)$. There is thus a trade-off between optimizing the $Cp^{-m/2}$ term in Theorem 2a and optimizing the $(C')^p$ term; the optimal choice of $\gamma$ depends upon the other inputs.

Notice, also, that for fixed $\gamma$ and for values of $q^\vee_i$ bounded away from zero, the constant $C$ is essentially a constant multiple of the bound on ${\rm conc}(AX)$ in Theorem 1. In fact, for (say) $\gamma=1$, we have
$$C\leq (.657)^m \prod_{i=1}^m \frac{1-q^\vee_i}{\sqrt{q^\vee_i}},$$
suggesting that the results of Theorem 2a are significantly better than those of Theorem 1 when $p$ is large enough that the $Cp^{-m/2}$ term dominates.

{\sl Proof of Theorem 2b.}\quad 
Recall that
\begin{align*}
c_i = \frac{1}{\gamma^2}\ln\Big[1+\alpha^\vee_i \big(1-\cos\frac{\gamma}{\sqrt{\alpha^\vee_i}}\big)\Big] \quad\quad &{\rm if}~\alpha^\vee_i\geq \frac{\gamma^2}{\pi^2}, \\
c_i = \frac{1}{\alpha^\vee_i\pi^2}\ln\Big[1+2\alpha^\vee_i\Big] \quad\quad &{\rm if}~\alpha^\vee_i\leq \frac{\gamma^2}{\pi^2},
\end{align*}
$$C := \prod\limits_{i=1}^m (2\pi c_i\alpha^\vee_i)^{-1/2},$$
and
$$C' := \max\limits_{1\leq i\leq m} e^{-\gamma^2 c_i/2}.$$

Regarding $c_i$ as a function of $\alpha^\vee_i$, we claim that this function is minimized at $\alpha^\vee_i=\frac{\gamma^2}{\pi^2}$. To demonstrate this claim, it suffices to check that:
\begin{enumerate}
\item The function $f(x) := \frac{\ln(1+2x)}{x}$ is decreasing for $0<x\leq\frac{\gamma^2}{\pi^2}$.
\item The function $g(x) := x(1-\cos\frac{\gamma}{\sqrt{x}})$ is increasing for $\frac{\gamma^2}{\pi^2}\leq x<\infty$.
\end{enumerate}

{\sl Proof of (1):} Differentiating, we obtain $f'(x) = \frac{1}{x^2} \big[\frac{2x}{1+2x} - \ln(1+2x)]$. In general, $\ln(1+u) > \frac{u}{1+u}$ for $u>0$, so we have $f'(x)<0$ for all $x>0$. In particular, $f(x)$ is decreasing for $0<x\leq\frac{\gamma^2}{\pi^2}$.

{\sl Proof of (2):} Differentiating, we obtain $g'(x) = 1 - \cos\frac{\gamma}{\sqrt{x}} - \frac{\gamma}{2\sqrt{x}} \sin \frac{\gamma}{\sqrt{x}}$. It will be convenient to define $y:=y(x)=\frac{\gamma}{\sqrt{x}}$. This change of variable bijectively transforms the interval $\frac{\gamma^2}{\pi^2}\leq x<\infty$ into the interval $0<y\leq\pi$. We may hence write $g'(x) = h(y)$, where
$$h(y) := 1 - \cos y - \frac{y}{2}\sin y.$$
Differentiating twice {\sl with respect to $y$}, we obtain
$$\frac{dh}{dy} = \frac{1}{2}\sin y - \frac{y}{2}\cos y \quad{\rm and}\quad
\frac{d^2 h}{dy^2} = \frac{y}{2}\sin y.$$
In particular, note that $h(0)=0$, $h'(0)=0$, and $h''(y)>0$ for $0<y<\pi$. It follows that $h(y)>0$ for $0<y\leq\pi$. Equivalently, $g'(x)>0$ (and $g(x)$ is increasing) for $\frac{\gamma^2}{\pi^2}\leq x<\infty$.

We have thus proved that $c_i$ is minimized when $\alpha^\vee_i = \frac{\gamma^2}{\pi^2}$, in which case\break
$c_i = \frac{1}{\gamma^2}\ln\big(1+\frac{2\gamma^2}{\pi^2} \big)$. That is to say,
$$c_i \geq \frac{1}{\gamma^2}\ln\big(1+\frac{2\gamma^2}{\pi^2} \big)$$
for all values of $\alpha^\vee_i$. It follows that
\begin{align*}
C \quad=\quad \prod\limits_{i=1}^m (2\pi c_i\alpha^\vee_i)^{-1/2}
\quad&\leq\quad \prod\limits_{i=1}^m \Big(\frac{2\pi}{\gamma^2} \ln\big(1+\frac{2\gamma^2}{\pi^2} \big) \cdot \frac{2q^\vee_i}{(1-q^\vee_i)^2} \Big)^{-1/2} \\
\quad&=\quad \Bigg[\frac{\gamma}{2\sqrt{\pi\ln\big(1+\frac{2\gamma^2}{\pi^2}\big)}} \Bigg]^m \prod_{i=1}^m \frac{1-q^\vee_i}{\sqrt{q^\vee_i}}
\end{align*}
and
\begin{align*}
C' \quad=\quad \max\limits_{1\leq i\leq m} e^{-\gamma^2 c_i/2} \quad&\leq\quad \exp\bigg(-\frac{\ln\big(1+\frac{2\gamma^2}{\pi^2}\big)}{2}\bigg) \\
\quad&=\quad \frac{1}{\sqrt{1+\frac{2\gamma^2}{\pi^2}}},
\end{align*}
proving Theorem 2b. $\blacksquare$

\section{Proof of Theorem 3}

We obtain Theorem 3 as a corollary of Proposition 3a, a more general result to follow. In order to state and prove Proposition 3a, we borrow the following notions from the theory of partially ordered sets (posets).
\begin{Defs}
Let $S$ be a poset and $x,y\in S$. We say that $x$ {\rm covers} $y$ if $x>y$ and $x\geq z\geq y\Rightarrow z\in\{x,y\}$.

A {\rm rank function} on a finite poset $S$ is a function ${\rm rk}:S\rightarrow {\Bbb Z}_{\geq 0}$, such that for all $x,y\in S$, if $x$ covers $y$, then ${\rm rk}(x)={\rm rk}(y)+1$. We say that ${\rm rk}(x)$ is the rank of element $x$. A {\rm layer} of a ranked poset is a level set of the rank function.

The chain of cardinality $N$ is denoted by $[N]$, and is automatically assigned herein the unique rank function which assigns its least element rank $0$. The product of two ranked posets $S,S'$ is automatically assigned rank function equal to the sum of the rank functions of $S,S'$.

An {\rm antichain} in a poset is a collection of pairwise incomparable elements. The {\rm width} of a poset $S$, denoted by $w(S)$, is the cardinality of its largest antichain(s). The {\rm Whitney number} $W_i$ of a ranked poset is the cardinality of its layer of rank $i$. If the width of a ranked poset is equal to its largest Whitney number, then we say that the poset has the {\rm Sperner property}.
\end{Defs}

For example, the ``Boolean cube'' $([2]\times[2]\times[2])$ has Whitney numbers $1,3,3,1$ and width $3$. Note that the width of any poset is greater than or equal to its largest Whitney number, because all layers are necessarily antichains.

Now we are ready to state
\begin{Prop3a}
Let $X_1,X_2,\ldots,X_p$ be independent, integer-valued random variables such that
$${\rm conc}(X_j)\leq {1\over N_j} \quad\quad {\rm for}~1\leq j\leq p,$$
where $N_1,N_2,\ldots,N_p$ are positive integers. Then
$${\rm conc}(X_1+\cdots+X_p) \leq {w\big([N_1]\times\cdots\times [N_p]\big)\over N_1N_2\cdots N_p}.$$
Moreover, given any fixed $N$ such that $2\leq N_1,N_2,\ldots,N_p<N$, we have
$${w\big([N_1]\times\cdots\times [N_p]\big)\over N_1N_2\cdots N_p} \sim \Big(\frac{\pi}{6} \sum\limits_{j=1}^p (N_j^2-1)\Big )^{-1/2}$$
as $p\rightarrow\infty$.
\end{Prop3a}

This proposition will be easiest to prove under the assumption that each $X_j$ is uniformly supported on $N_j$ points (with mass $\frac{1}{N_j}$ at each). To justify passing to this case, we will use the following definition, and the two lemmas after it:

\begin{Def}
A discrete random variable $Y$ is a {\rm mixture} of random variables $Y_1,Y_2,\ldots$ if its probability mass function lies in the convex hull of the probability mass functions of $Y_1,Y_2,\ldots$.
\end{Def}

\begin{Lem}
Let $Y$ be a random variable, supported on ${\Bbb Z}_{\geq 0}$, such that ${\rm conc}(Y)\leq \frac{1}{N}$. Then $Y$ can be written as a mixture of random variables $Y_1,Y_2,\ldots$, such that each $Y_k$ is uniformly supported on $N$ points, i.e., has an $N$-point support with probability mass $\frac{1}{N}$ at each point in its support.
\end{Lem}

{\sl Proof of Lemma 9.}\quad
Let ${\mathcal M}$ be the space of probability measures on ${\Bbb Z}_{\geq 0}$. Let
$${\mathcal M}(N) := \bigg\{\mu\in{\mathcal M}:~~\max_k \mu(\{k\}) \leq \frac{1}{N}\bigg\}$$
and
$${\mathcal M}_{\rm u}(N) := \{\mu\in{\mathcal M}:~~\mu~{\rm is~uniformly~supported~on}~N~{\rm points}\}.$$
By the Krein-Milman theorem, ${\mathcal M}(N)$ is the convex hull of its extreme points. We claim that the extreme points are precisely the points of ${\mathcal M}_{\rm u}(N)$. It is immediately evident that each point of ${\mathcal M}_{\rm u}(N)$ is an extreme point of ${\mathcal M}(N)$. Conversely, suppose $\mu\in{\mathcal M}(N)\backslash{\mathcal M}_{\rm u}(N)$. Thus there is some $k\in{\Bbb Z}_{\geq 0}$ such that $0<\mu(\{k\})<\frac{1}{N}$, but in fact, there must be at least two distinct such $k$, since the total mass of $\mu$ is 1 (an integer multiple of $\frac{1}{N}$). Therefore, $\mu$ is not an extreme point of ${\mathcal M}(N)$.

This proves our claim. Hence the probability measure associated to $Y$ can be written as a countable convex combination of points of ${\mathcal M}_{\rm u}(N)$, each of which defines the distribution of a random variable $Y_k$ (proving the lemma). $\square$

\begin{Lem}[Properties of superpositions]
If $Y$ is a mixture of random variables $Y_1,Y_2,\ldots$, then:
\begin{enumerate}
\item There is some $k\geq 1$ for which ${\rm conc}(Y) \leq {\rm conc}(Y_k)$.
\item If $Z$ is a random variable and $f$ a function such that $Z=f(Y)$, then $Z$ is a mixture of random variables $Z_1,Z_2,\ldots$, where $Z_k = f(Y_k)$.
\end{enumerate}
\end{Lem}

{\sl Proof of Lemma 10.}\quad
By the definition of {\sl mixture}, there exist nonnegative $\alpha_1,\alpha_2,\ldots$ such that $\alpha_1+\alpha_2+\cdots=1$ and such that
$${\bf Pr}[Y=y] = \sum_{k=1}^\infty \alpha_k{\bf Pr}[Y_k=y].$$
Thus by the pigeonhole principle, for arbitrary $y$, there exists $k=k(y)$ such that
$${\bf Pr}[Y=y] \leq {\bf Pr}[Y_k=y].$$
Choosing $y$ such that ${\rm conc}(Y) = {\bf Pr}[Y=y]$, we conclude that ${\rm conc}(Y) \leq {\rm conc}(Y_k)$ for this $k$. This proves claim (1) in the lemma. Claim (2) is self-evident. $\square$

The heart of the proof of Proposition 3a is the following version of the local limit theorem:

\begin{Def}
A sequence $(\ldots,b_{-1},b_0,b_1,b_2,\ldots)$ of nonnegative real numbers is {\sl properly log-concave} if it is log-concave (i.e., $b_{t-1}b_{t+1}\leq b_t^2$ for all $t$) and has no internal zeroes (i.e., if $b_t>0$ and $b_{t+k}>0$, then $b_{t+1},b_{t+2},\ldots,b_{t+k-1}>0$).
\end{Def}

\begin{Lem}[Bender]
Suppose that $\big(\zeta_p: p\in{\Bbb N}\big)$ is a sequence of integer-valued random variables, $\big(F_p\big)$ are the corresponding distribution functions, and $\big(\sigma_p\big)$ and $\big(\mu_p\big)$ are sequences of real numbers such that $\lim\limits_{p\rightarrow\infty} F_p(\sigma_p x+\mu_p) = \frac{1}{\sqrt{2\pi}} \int_{-\infty}^x e^{-t^2/2} dt$ for every $x\in{\Bbb R}$. Also suppose that $\sigma_p\rightarrow\infty$ as $p\rightarrow\infty$. Further, suppose that, for every $p$, the sequence $b_p(t):={\bf Pr}(\zeta_p=t)$ is properly log-concave with respect to $t$. Then
$$\lim_{p\rightarrow\infty} \sigma_p{\bf Pr}\big(\zeta_p = \lfloor\sigma_p x+\mu_p\rfloor\big) = \frac{1}{\sqrt{2\pi}} e^{-x^2/2}$$
uniformly for all $x\in{\Bbb R}$.
\end{Lem}

This result originally appeared in~\cite{Bender1973}, but the above statement is based on its treatment in~\cite{Engel1997}; see either source for a proof.

{\sl Proof of Proposition 3a.}
For $j=1,2,\ldots,p$, we are given to assume that\break
${\rm conc}(X_j)\leq\frac{1}{N_j}$. By Lemma 9, each $X_j$ is a superposition of some random variables which are each uniformly supported on some $N_j$ points. Thus the random vector $X=(X_1,\ldots,X_p)$ is a mixture of random vectors each of the form $X^{(k)}:=(X^{(k)}_1,\ldots,X^{(k)}_p)$, where the coordinates are independent and each $X^{(k)}_j$ is uniformly supported on $N_j$ points. The sum $X_1+\cdots+X_p$ is a function of $X$, so by using both parts of Lemma 10, we see that
$${\rm conc}(X_1+\cdots+X_p) \leq {\rm conc}(X^{(k)}_1+\cdots+X^{(k)}_p)$$
for some $k$. Since we are seeking an upper bound on ${\rm conc}(X_1+\cdots+X_p)$, we assume with no loss of generality that $X=X^{(k)}$, or, more to the point, that each coordinate $X_j$ is uniformly supported on $N_j$ points (with mass $\frac{1}{N_j}$ on each).

Denote the support of $X_j$ by $\{a_{j1},a_{j2},\ldots,a_{j{N_j}}\}$, where $a_{j1} < a_{j2} <\cdots< a_{j{N_j}}$. Then
$$a_{1{i_1}} + a_{2{i_2}} + \cdots + a_{p{i_p}} = a_{1{i'_1}} + a_{2{i'_2}} + \cdots + a_{p{i'_p}}$$
implies that the $p$-tuples $(i_1,i_2,\ldots,i_p)$ and $(i'_1,i'_2,\ldots,i'_p)$ are identical or incomparable in
$[N_1]\times\cdots\times [N_p]$. It follows that
$${\rm conc}(X_1+\cdots+X_p) \leq {w\big([N_1]\times\cdots\times [N_p]\big)\over N_1N_2\cdots N_p}.$$
This proves the first claim of Proposition 3a.

For the remainder of the proof, assume that $2\leq N_1,N_2,\ldots,N_p<N$ for some integer $N$. We are going to apply Lemma 13. Let $\zeta_p$ denote the rank of a uniformly distributed random element of $[N_1]\times[N_2]\times\cdots\times[N_p]$. Set $\mu_p:=\frac{N_1+\cdots+N_p}{2}$ and $\sigma_p^2=\sum_{j=1}^p \frac{N_j^2-1}{12}$. It is easily verified that $\mu_p$ and $\sigma_p^2$ are respectively the mean and the variance of $\zeta_p$. By Lyapunov's central limit theorem~\cite{Billingsley}, the condition
$$\lim\limits_{p\rightarrow\infty} F_p(\sigma_p x+\mu_p) = \frac{1}{\sqrt{2\pi}} \int_{-\infty}^x e^{-t^2/2} dt$$
in Lemma 11 is satisfied. The hypothesis $\sigma_p\rightarrow\infty$ is plainly also satisfied.

To see that the sequence $b_p(t):={\bf Pr}(\zeta_p=t)$ is properly log-concave, we note that this sequence is proportional to the Whitney numbers of the chain product $[N_1]\times[N_2]\times\cdots\times[N_p]$, which is the convolution of the sequences of Whitney numbers for the factor chains. Each factor chain has Whitney numbers $1,1,\ldots,1,0,0,\ldots$ (a properly log-concave sequence). Furthermore, the convolution of properly log-concave sequences is again properly log-concave, see e.g.~\cite{Kook2006}. Thus, $\big(b_p(t)\big)$ is properly log-concave.

All conditions of Lemma 11 have been verified, so the conclusion holds:
$$\lim_{p\rightarrow\infty} \sigma_p{\bf Pr}\big(\zeta_p = \lfloor\sigma_p x+\mu_p\rfloor\big) = \frac{1}{\sqrt{2\pi}} e^{-x^2/2}$$
uniformly for all $x\in{\Bbb R}$. Setting $x=0$, we obtain
\begin{align*}
{\bf Pr}\big(\zeta_p=\lfloor\mu_p\rfloor\big) &\sim \frac{1}{\sqrt{2\pi\sigma_p}} \\
&= \Big(\frac{\pi}{6} \sum\limits_{j=1}^p (N_j^2-1)\Big )^{-1/2}.
\end{align*}

Finally, we observe that chain products have the Sperner property~\cite{Engel1997}. In particular, the width in the above formula is equal to the Whitney number $W_{\lfloor \mu_p \rfloor}$, so that
$${w\big([N_1]\times\cdots\times [N_p]\big)\over N_1N_2\cdots N_p} = {\bf Pr}\big(\zeta_p=\lfloor\mu_p\rfloor\big).$$
This completes the proof of the proposition. $\blacksquare$

As an instance of Proposition 3a, we derive Theorem 3:

{\sl Proof of Theorem 3.}\quad
As noted in the proof of Theorem 1, we have ${\rm conc}(X_j{\bf a}_j)=\frac{1}{{\bf E}(X_j)+1} \leq \frac{1}{\lfloor{\bf E}(X_j)+1\rfloor}$ for $1\leq j\leq n$. Since ${\bf a}_1,{\bf a}_2,\ldots,{\bf a}_m$ are linearly independent, we have
\begin{align*}
{\rm conc}(AX) &= \prod_{i=1}^m {\rm conc}(X_i{\bf a}_i + X_{m+i}{\bf a}_i + X_{2m+i}{\bf a}_i \cdots + X_{(p-1)m+i}{\bf a}_i) \\
&= \prod_{i=1}^m {\rm conc}(X_i+X_{m+i}+X_{2m+i}+\cdots+X_{(p-1)m+i}) \\
&\lesssim \prod_{i=1}^m \Big( \frac{\pi p}{6}\big(\lfloor{\bf E}(X_i)+1\rfloor^2 -1\big)\Big)^{-1/2},
\end{align*}
where the last claim follows by Proposition 3a. Finally, by section 2.1, (2), we infer Theorem 3. $\blacksquare$

\section{Acknowledgments.}
The author thanks Alexander Barvinok and Roman Vershynin for fruitful discussions, and the former also for making many helpful comments on the style of this paper.

\bibliographystyle{plain}
\bibliography{concentration_bibliography}

\begin{thebibliography}{10}

\bibitem{BH2010}
A.~Barvinok and J.~Hartigan.
\newblock {Maximum entropy Gaussian approximation for the number of integer
  points and volumes of polytopes}.
\newblock {\em Advances in Applied Mathematics}, 45:252--289, 2010.

\bibitem{Bender1973}
E.~A. Bender.
\newblock {Central and local limit theorems applied to asymptotic enumeration}.
\newblock {\em Journal of Combinatorial Theory, Series A}, 15:91--111, 1973.

\bibitem{Billingsley}
P.~Billingsley.
\newblock {\em {Probability and Measure, 2nd ed.}}
\newblock Wiley, New York, 1986.

\bibitem{Burchard}
A.~Burchard.
\newblock {A short course on rearrangement inequalities}.
\newblock \url{http://www.math.utoronto.ca/almut/rearrange.pdf}.

\bibitem{DeLoera2005}
J.~De~Loera.
\newblock The many aspects of counting lattice points in polytopes.
\newblock {\em Mathematische Semesterberichte}, 52:175--195, 2005.

\bibitem{DeLoera2009}
J.~De~Loera.
\newblock {Counting and estimating lattice points: tools from algebra,
  analysis, convexity, and probability}.
\newblock {\em Optima: Newsletter of the Mathematical Programming Society},
  Dec. 2009.

\bibitem{DE1985}
P.~Diaconis and B.~Efron.
\newblock Testing for independence in a two-way table: new interpretations of
  the chi-square statistic.
\newblock {\em Annals of Statistics}, 13:845--874, 1985.

\bibitem{Engel1997}
K.~Engel.
\newblock {\em {Sperner Theory}}.
\newblock Cambridge University Press, Cambridge, 1997.

\bibitem{Erdos1945}
P.~Erd\"os.
\newblock {On a lemma of Littlewood and Offord}.
\newblock {\em Bulletin of the American Mathematical Society}, 51:898--902,
  1945.

\bibitem{GJ1979}
M.R. Garey and S.J. Johnson.
\newblock {\em {Computers and Intractability: A Guide to the Theory of
  NP-Completeness}}.
\newblock Freeman, San Francisco, 1979.

\bibitem{Halasz1977}
G.~Hal\'asz.
\newblock {Estimates for the concentration function of combinatorial number
  theory and probability}.
\newblock {\em Periodica Mathematica Hungarica}, 8:197--211, 1977.

\bibitem{Howard}
R.~Howard.
\newblock {Estimates on the concentration function of sets in ${\Bbb R}^d$:
  Notes on lectures of Oskolkov}.
\newblock \url{http://www.math.sc.edu/~howard/Notes/concentration.pdf}.

\bibitem{JSV2004}
M.~Jerrum, A.~Sinclair, and E.~Vigoda.
\newblock A polynomial-time approximation algorithm for the permanent of a
  matrix with nonnegative entries.
\newblock {\em Journal of the ACM}, 51:671--697, July 2004.

\bibitem{Kook2006}
W.~Kook.
\newblock {On the product of log-concave polynomials}.
\newblock {\em INTEGERS: Electronic Journal of Combinatorial Number Theory}, 6,
  2006.

\bibitem{RV2010}
M.~Rudelson and R.~Vershynin.
\newblock {The Littlewood-Offord problem and invertibility of random matrices}.
\newblock {\em Advances in Mathematics}, 218:600--633, 2008.

\bibitem{TV2009}
T.~Tao and V.~Vu.
\newblock Inverse littlewood-offord theorems and the condition number of random
  discrete matrices.
\newblock {\em Annals of Mathematics}, 169:595--632, 2009.

\end{thebibliography}

\end{document}